\documentclass[letterpaper, 10 pt, conference]{ieeeconf}  
\usepackage{graphicx}
\usepackage{amsmath}
\usepackage{amssymb}
\usepackage{CJK}
\usepackage{color}
\usepackage{dsfont}
\usepackage{subfigure}
\usepackage{booktabs}
\usepackage{microtype}
\usepackage{cite,url}

\usepackage{graphicx}
\usepackage{amssymb,amsfonts,eucal}
\usepackage{cite,url}
\usepackage{mathrsfs}
\usepackage{bm}
\usepackage{float}
\usepackage{amsmath}
\usepackage{subfigure}
\usepackage{pstool}
\usepackage{color}
\usepackage{array}
\usepackage{cancel}

\allowdisplaybreaks

\title{\LARGE \bf Stochastic Optimal Power Flow Based on Data-Driven Distributionally Robust Optimization}

\author{Yi Guo$^\dag$,
        Kyri Baker$^\ddag$,
        Emiliano Dall'Anese$^*$, 
        Zechun Hu$^{**}$,
        Tyler Summers$^\dag$ ~
\thanks{$^\dag$Y. Guo and T. Summers are with the Department
of Mechanical Engineering, The University of Texas at Dallas, Richardson,
TX, USA, email: \{yi.guo2,tyler.summers\}@utdallas.edu. $^\ddag$K. Baker is with the Department of Civil, Envirnomental, and Architectural Engineering, The University of Colorado, Boulder, CO, USA, email:kyri.baker@colorado.edu. $^*$E. Dall'Anese is with the National Renewable Energy Laboratory, Golden, CO, USA, email: emiliano.dallanese@nrel.gov. $^{**}$Z. Hu is with the Department of Electrical Engineering, Tsinghua University, Beijing, P.R. China, email: zechhu@tsinghua.edu.cn. This research is supported by NSF under grant CNS-1566127.}}
\IEEEoverridecommandlockouts
\begin{document}
\maketitle
\thispagestyle{empty}
\pagestyle{empty}
\maketitle

\begin{abstract}
We propose a data-driven method to solve a stochastic optimal power flow (OPF) problem based on limited information about forecast error distributions. The objective is to determine power schedules for controllable devices in a power network to balance operation cost and conditional value-at-risk (CVaR) of device and network constraint violations. These decisions include scheduled power output adjustments and reserve policies, which specify planned reactions to forecast errors in order to accommodate fluctuating renewable energy sources. Instead of assuming the uncertainties across the networks follow prescribed probability distributions, we assume the distributions are only observable through a finite training dataset. By utilizing the Wasserstein metric to quantify differences between the empirical data-based distribution and the real data-generating distribution, we formulate a distributionally robust optimization OPF problem to search for power schedules and reserve policies that are robust to sampling errors inherent in the dataset. A multi-stage closed-loop control strategy based on model predictive control (MPC) is also discussed. A simple numerical example illustrates inherent tradeoffs between operation cost and risk of constraint violation, and we show how our proposed method offers a data-driven framework to balance these objectives.
\end{abstract}


%
\IEEEpeerreviewmaketitle

\section{Introduction}
The continued integration of renewable energy sources (RESs) in power systems is making it more complicated for system operators to balance economic efficiency and system reliability and security. As penetration levels of RESs reach substantial fractions of total supplied power, networks will face high operation risks under current operational paradigms. As it becomes more difficult to predict the net load, large forecast errors can lead to power quality and reliability issues causing significant damage and costly outages. Future power networks will require more sophisticated methods for managing these risks, at both transmission and distribution levels.

The flexibility of controllable devices, including power-electronics-interfaced RESs, can be utilized to balance efficiency and risk with  optimal power flow methods \cite{opf,dommel,alsac,Baldick,low1,low2}, which aim to determine power schedules for controllable devices in a power network to optimize an objective function. However, most OPF methods in the research literature and those widely used in practice are deterministic, assuming point forecasts of exogenous power injections and ignoring forecast errors. Increasing forecast errors push the underlying distributed feedback controllers that must handle the transients caused by these errors closer to stability limits \cite{kundur}.

More recently, research focus has turned to stochastic and robust optimal power flow methods that explicitly incorporate forecast errors, in order to more systematically trade off efficiency and risk and to ease the burden on feedback controllers \cite{yong,capitanescu,conejo,bienstock,Vrakopoulou1,vrakopoulou2,warrington,zhang1,perninge,roald1,tyler,jabr,roald2,roald3,LiMathieu,zhang2,lubin,Baker,DallAnese1}. Many formulations assume that uncertain forecast errors follow a prescribed probability distribution (commonly, Gaussian \cite{bienstock,lubin,roald1}) and utilize analytically tractable reformulations of probabilistic constraints. In practice, forecast error probability distributions are never known; they are only observed indirectly through finite datasets. Sampling-based methods have been applied with a focus on quantifying the probability of constraint violation \cite{Vrakopoulou1,vrakopoulou2} and for constraining or optimizing conditional value at risk (CVaR) \cite{zhang1,tyler,DallAnese1}. Distributionally robust approaches use data to estimate distribution parameters (e.g., mean and variance) and aim to be robust to any data-generating distribution consistent with these parameters \cite{tyler,roald3,LiMathieu,Baker,DallAnese1}. Others take a robust approach, assuming only knowledge of bounds on forecast errors and enforcing constraints for any possible realization, e.g., \cite{warrington,jabr}. Overall, this line of recent research has explored tractable approximations and reformulations of difficult stochastic optimal power flow problems. However, none of the existing work explicitly accounts for sampling errors arising from limited data, which in operation can cause poor out-of-sample performance\footnote{Out-of-sample performance is an evaluation of the optimal decisions using a dataset that is different from the one used to obtain the decision, which can be tested with Monte Carlo simulation}. Even with sophisticated recent stochastic programming techniques, decisions can be overly dependent on small amounts of relevant data, a phenomenon akin to overfitting in statistical models.

We propose a data-driven method to solve a stochastic optimal power flow problem based on limited information about forecast error distributions available through a finite training dataset. Our approach is inspired by recent results in the optimization literature on data-driven distributionally robust optimization \cite{mohajerin}. We formulate a tractable data-driven distributionally robust optimal power flow problem where the Wasserstein metric \cite{givens} is used to quantify differences between an empirical data-based distribution and the real data-generating distribution. In contrast to previous work, we obtain power schedules and reserve  policies that are explicitly robust to sampling errors inherent in the dataset. This approach achieves superior out-of-sample performance guarantees in comparison to other stochastic optimization approaches, effectively regularizing against overfitting the decisions to limited available data. A numerical example based on a modified IEEE 118 bus test network illustrates inherent tradeoffs between operation cost and risk of constraint violation, and we show how our proposed method offers a data-driven framework to systematically balance these objectives The results also point out that making inaccurate assumptions about probability distributions can cause underestimation of risk.

The rest of the paper is organized as follows: Section II describes the network model and the stochastic OPF problem. Section III describes the data-driven distributionally robustness optimization framework and formulates the data-driven stochastic OPF problem. Section IV provides simulation results and discussion. Section V concludes the paper.

\textbf{Notation}: The inner product of two vectors $a,b \in \mathbb{R}^m$  is denoted by $\langle a, b \rangle := a^\textbf{T}b$. The $N$-fold product of distribution $\mathds{P}$ on a set $\Xi$ is denoted by $\mathds{P}^N$, which represents a distribution on the Cartesian product space $\Xi_{N} = \Xi \times ... \times \Xi$. Superscript `` $\hat{.}$ " is reserved for the objects that depend on a training dataset $\hat{\Xi}_N$. We use $(.)'$ to denote  vector or matrix transpose and $A^T$ to denote matrix $A$ to the power $T$.

\section{Network Model and Optimal Power Flow}
We consider $N_d$ devices connected via a transmission network within a planning time horizon \emph{T} discrete time steps. The elements connected to the transmission system may include 1) traditional generators; 2) fixed, deferrable, and curtailable loads; 3) storage devices like batteries and plug-in electric vehicles, which are able to act as both generators and loads. In this case, we model two types of devices: devices (e.g., conventional thermal and hydro generators, deferrable/curtailable loads and storage devices) with controllable power flow affected by decision variables; and devices (e.g., flexible and intermittent renewable energy or fixed loads) with fixed or uncertain power flow which will not be affected by decision variables.

\subsection{Devices with Controllable Power Flow}
The power flow of each controllable device is modeled with a discrete-time linear dynamical system
\begin{equation}\label{dynamic}
x_{t+1}^j = \bar{A}_j x_{t}^j + \bar{B}_j u_{t}^j,
\end{equation}
where device $j$ at time $t$ has internal state $x_t^j \in \mathbb{R}^{n_j}$, dynamics matrix $\bar{A}_j \in \mathbb{R}^{n_j \times n_j}$, input matrix $\bar{B}_j \in \mathbb{R} ^{n_j \times m_j}$, and control input $u_t^j \in \mathbb{R}^{m_j}$. The first element of $x_t^{j}$ corresponds to the power injection of device $j$ at time $t$ into the network at its bus, and other elements describe internal dynamics such as state-of-charge (SOC) of  energy storage devices. States and inputs over the planning horizon could be expressed in the following notation: $ \textbf{x}^j = [(x_1^j)', ..., (x_{T}^j)' ]' \in \mathbb{R}^{n_jT} $ and $\textbf{u}^j = [(u_0^j)', ..., (u_{T-1}^j)' ]' \in \mathbb{R}^{m_jT}$. We can then write compactly:
\begin{equation}\label{statefunction}
\textbf{x}^j = A_jx_0^{j} + B_j\textbf{u}^j,
\end{equation}
where
\begin{equation*}
A_j :=
\begin{bmatrix}
\bar{A}_j\\
\bar{A}_j^2\\
\vdots\\
\bar{A}_j^{T}
\end{bmatrix}
,
B_j :=
\begin{bmatrix}
\bar{B}_j & 0 & \dots & 0\\
\bar{A}_j\bar{B}_j & \bar{B_j}& \ddots & 0 \\
\vdots & \ddots & \ddots & \vdots \\
\bar{A}_j^{T-1} \bar{B}_j & \dots & \bar{A}_j\bar{B}_j & \bar{B}_j
\end{bmatrix}.
\end{equation*}

\subsection{Devices with Uncontrollable Power Flow}
The fixed power flow for device $j$ is given by $r_j + G_j\xi$ with positive values denoting net power injection into the network. The vector $r_j \in \mathbb{R} ^T$ denotes the predicted power injection over the planning horizon, and the matrix $G_j \in \mathbb{R}^{T \times N_\xi T}$ maps the random vector $\xi \in \Xi \subseteq \mathbb{R}^{N_\xi T}$ to a prediction error of the power injection or extraction for device $j$. If uncertainty of device $j$ is not explicitly modeled, then $G_j = 0$.

We will assume that the probability distribution $\mathbb{P}$ of $\xi$ is unknown and observed only through a finite dataset $\hat \Xi_N$.

\subsection{Cost Functions and Constraints}
An operating objective of our problem is to minimize a sum of cost functions associated with controllable devices $J_j : \mathbb{R}^{n_jT } \times \mathbb{R}^{m_jT} \to \mathbb{R}$. The cost functions are modeled as convex quadratic:
\begin{equation}\label{costfunction}
J_j(\textbf{x}^j, \textbf{u}^j) := f_{jx}'\textbf{x}^j +  \frac{1}{2}{\textbf{x}^j}'H_{jx}\textbf{x}^j + f_{ju}'\textbf{u}^j +  \frac{1}{2}{\textbf{u}^j}'H_{ju}\textbf{u}^j + c_j,
\end{equation}
wth  $H_{jx}$ and $H_{ju}$ being positive semidefinite matrices.\\
We consider: 1) power balance constraints, 2) power line flow constraints and 3) local device constraints. The power balance constraints are
\begin{equation}\label{powerbalance}
\sum_{j=1}^{N_d} (r_j + G_j\xi +C_j\textbf{x}^j) = 0,
\end{equation}
which ensures that the generation and consumption of power inside the transmission system are balanced over the planning horizon. The matrix $C_j$ selects the first element of $\textbf{x}^j$. In this paper, we employ a widely used linearization of the nonlinear AC power flow equations, in which it is assumed that lines are lossless, normalized voltage magnitudes are close to unity, and voltage angle differences are small \cite{tyler,warrington}. In this case, the line flow are linear functions of nodal injections.

The $L$ bi-directional transmission lines power flow constraints can be expressed as
\begin{equation} \label{lineconstraints}
\sum_{j=1}^{N_d} \Gamma_j(r_j + G_j\xi + C_j\textbf{x}^j) \le \bar{p},
\end{equation}
where $\Gamma_j \in \mathbb{R}^{2LT \times T}$ maps the power injection or extraction of each device to its contribution to each connected line and can be constructed from network line impedances, and $\bar{p}$ denotes nominal line flow limits. The local constraints are also modeled as linear inequalities of the form:
\begin{equation}\label{localdeviceconstraints}
T_j\textbf{x}^j + U_j\textbf{u}^j +Z_j\xi \leq w_j.
\end{equation}
where $T_j \in \mathbb{R}^{l_j \times n_jT}$, $U_j \in \mathbb{R}^{l_j \times m_jT}$, and $Z_j \in \mathbb{R}^{l_j \times N_\xi T}$, and $w_j\in \mathbb{R}^{l_j}$ is a local constraint parameter vector. These constraints can be used to model allowable power injection ranges and other device limits.

\subsection{Reserve Policies}
Deterministic OPF formulations ignore the prediction error $\xi$ and compute an open-loop input sequence for each device. In a stochastic setting, one must optimize over causal \emph{policies} $\textbf{u}^j = \pi_j(\xi)$, where $\pi_j : \mathbb{R}^{N_\delta T} \to \mathbb{R}^{m_jT}$ is a measurable function that specifies how each device should respond to forecast errors as they are discovered. We can now formulate a finite horizon stochastic optimization problem
\begin{equation} \label{stochopf}
\begin{split}
& J_{Cost} =  \inf_{\pi_j \in \Pi_c}  \mathds{E}\sum_{j=1}^{N_d} J_j(A_jx_0^{j} + B_j\pi_j(\xi), \pi_j(\xi)),\\
& \text{subject to}\\
  &\sum_{j=1}^{N_d} (r_j + G_j\xi + C_j(A_jx_0^{j} + B_j\pi_j(\xi))) = 0, \forall\xi,\\
  & \mathds{E}f_1\bigg(\sum_{j=1}^{N_d} \Gamma_j(r_j + G_j\xi + C_j (A_jx_0^j + B_j\pi_j(\xi))) -\bar{p} \bigg) \le 0 ,\\
  & \mathds{E}f_2\bigg(T_j(A_jx_0^j + B_j\pi_j(\xi)) + U_j\pi_j(\xi) - w_j \bigg) \le 0, j = 1,... N_d.
\end{split}
\end{equation}
where $f_1$ and $f_2$ denote general constraint risk functions, which in our case will depend on forecast error data and auxiliary optimization variables. Also, the cost function is proportional to the first and second moments of the uncertainties $\xi$, because each device cost function are convex quadratic. The problem \eqref{stochopf} is infinite dimensional, so we restrict attention to affine policies:
\begin{equation} \label{inputpolicy}
\textbf{u}^j = D_j\xi + e_j.
\end{equation}
so that each participant device $j$ (e.g., traditional generators, flexible loads or energy storage) power schedule $\textbf{u}^j$ is parameterized by a nominal schedule $e_j = [e_0^j, ... , e_{T-1}^j]'$ plus a linear function $D_j$ of prediction error realizations. To obtain causal policies, $D_j$ must be lower-triangular. The $D_j$ matrices can be interpreted in terms of planned Automatic Generation Control (AGC) parameters \cite{warrington}. Under affine policies, the power balance constraints are linear functions of the distribution of $\xi$, which are equivalent to
\begin{equation} \label{constant equivalence}
\sum_{j=1}^{N_d}(r_j + C_j(A_jx_0^j + B_je_j)) = 0, \sum_{j=1}^{N_d}(G_j + C_jB_jD_j) = 0.
\end{equation}

% ...

\section{Data-Driven Distributionally Robust Optimization Via the Wasserstein Metric}
There is a variety of ways to reformulate the general stochastic OPF problem \eqref{stochopf} to obtain tractable problems that can be solved by standard convex optimization solvers. These include assuming specific functional forms for the forecast error distribution (e.g., Gaussian) and using specific constraint risk functions, such as those encoding value at risk (i.e., chance constraints), conditional value at risk, distributional robustness, and support robustness. In all cases, the out-of-sample performance of the resulting decisions in operational practice ultimately relies on the quality of data describing the forecast errors and the validity of assumptions made about probability distributions. Many existing approaches make either too strong or too weak assumptions that lead to underestimation or overestimation of risk. In this section, we review a recently proposed tractable method for data-driven distributionally robust optimization \cite{mohajerin}. We then use it to formulate a data-driven distributionally robust OPF problem that is based exclusively on a finite training dataset $\hat{\Xi}_N$, is robust to sampling errors, gives explicit control of conservativeness, and offers superior out-of-sample performance guarantees. This is accomplished by optimizing over the worst-case distribution within a Wasserstein ball in the space of probability distributions centered at the uniform distribution on the set of training data. This turns out to admit a tractable reformulation for certain risk functions, including CVaR.

\subsection{Data-Driven Stochastic Programming}
Consider the stochastic program
\begin{equation} \label{normalstochasticprogram}
J^* := \inf_{y \in \mathbb{Y}} \bigg\{\mathds{E}^{\mathds{P}}[h(y, \xi)] := \int_{\Xi} h(y, \xi)\mathds{P}(d\xi) \bigg\},
\end{equation}
with decision variable $y \in \mathbb{Y} \subseteq \mathbb{R}^n$, random vector $\xi$ with probability distribution $\mathds{P}$ supported on $\Xi \subseteq \mathbb{R}^m$, and cost function $h : \mathbb{R}^n \times \mathbb{R}^m \to \bar{\mathbb{R}}$. In virtually all practical settings, the distribution $\mathds{P}$ is unknown, so the problem is in some sense ill-posed, since it is not even possible to evaluate the cost function for a given decision. We typically have only a finite training dataset of $N$ independent samples observed from $\mathds{P}$, which we denote by $\hat{\Xi}_N := \{\hat{\xi_i}\}_{i \leq N} \subseteq \Xi$. 

The out-of-sample performance of a decision $\hat{y} \in \mathbb{Y}$ is $\mathds{E}^\mathds{P}[h(\hat{y}, \xi)]$, which can be interpreted as the expected cost of $\hat{y}$ under a new sample of $\xi$ that is independent of the training dataset. In any stochastic optimization problem, we seek a decision with good out-of-sample performance. A common approach is to approximate $\mathds{P}$ based on the training data $\hat \Xi_N$, effectively replacing $\mathds{P}$ with a discrete empirical distribution defined by the uniform distribution on $\hat \Xi_N$, which is denoted $\hat{\mathds{P}}_N$. However, this can cause decisions to be overly attuned to the data and lead to poor out-of-sample performance, especially when data is high-dimensional, limited in quantity, and expensive to obtain. The same is true for other approaches that utilize compressed summaries of $\hat \Xi_N$ (e.g., mean and variance) or assume a specific distribution that generated $\hat \Xi_N$.

Instead, we adopt an approach that explicitly accounts for our ignorance of true data-generating distribution $\mathds{P}$, is based explicitly on the available data $\hat \Xi_N$, and provides out-of-sample performance guarantees \cite{mohajerin}. A set 
$\hat{\mathcal{P}}_N$ of distributions is constructed from $\hat \Xi_N$ that contains those that could most plausibly have generated $\hat \Xi_N$. One can then formulate a distributionally robust optimization (DRO) problem:
\begin{equation}\label{generalDRO}
\hat{J}_{DRO} := \inf_{y\in \mathbb{Y}} \sup_{\mathds{Q}\in\hat{\mathcal{P}}_N} \mathds{E}^{\mathds{Q}}[h(y,\xi)],
\end{equation}
where the objective is to minimize the expected cost function with respect to the \emph{worst-case} data-generating distribution contained in $\hat{\mathcal{P}}_N$.

\subsection{The Wasserstein Metric}
The ambiguity set $\hat{\mathcal{P}}_N$ is defined using the Wasserstein metric, which defines a distance in the space $\mathcal{M}(\Xi)$ of all probability distributions $\mathds{Q}$ supported on $\Xi$ with $\mathds{E}^\mathds{Q}[||\xi||] = \int_\Xi ||\xi|| \mathds{Q}(d\xi) < \infty$.

\textbf{Definition} (Wasserstein Metric). Denote by $\mathcal{L}$ the space of all Lipschitz continuous functions $f: \Xi \to \mathbb{R}$ with Lipschitz constant less than or equal to 1. The Wasserstein metric $d_W : \mathcal{M}(\Xi) \times \mathcal{M}(\Xi) \to \mathbb{R}$ is defined as
\begin{equation} \label{WassersteinMetric}
\begin{array}{l}
d_W(\mathds{Q}_1, \mathds{Q}_2)
= \sup_{f \in \mathcal{L}}\bigg( \int_\Xi f(\xi)\mathds{Q}_1(d\xi) - \int_\Xi f(\xi)\mathds{Q}_2(d\xi) \bigg), \\
\forall \mathds{Q}_1, \mathds{Q}_2 \in \mathcal{M}(\Xi).
\end{array}
\end{equation}
Intuitively, the Wasserstein metric quantifies the minimum ``transportation'' cost to move mass from one distribution to another. We can now use the Wasserstein metric to define the ambiguity set
\begin{equation} \label{ambiguityset}
\hat{\mathcal{P}}_N := \bigg\{ \mathds{Q} \in \mathcal{M}(\Xi): d_w(\hat{\mathds{P}}_N, \mathds{Q}) \leq \varepsilon \bigg\},
\end{equation}
which contains all distributions within a Wasserstein ball of radius $\varepsilon$ centered at the empirical distribution $\hat{\mathds{P}}_N$. The radius $\varepsilon$ can be chosen so that the ball contains the true distribution $\mathds{P}$ with a prescribed confidence level and leads to performance guarantees \cite{mohajerin}. The radius $\varepsilon$ also explicitly controls the conservativeness of the resulting decision; large $\varepsilon$ will produce decisions that rely less on the specific features of the dataset $\hat \Xi_N$ and give better robustness to sampling errors.

\subsection{Distributionally Robust Optimization of CVaR}
The distributionally robust optimization problem \eqref{generalDRO} does not seem any easier than \eqref{normalstochasticprogram}, and seems quite possibly even worse. However, it was shown in \cite{mohajerin} that the \eqref{generalDRO} admits a tractable reformulation and in fact for piecewise affine cost functions of the form $h(y,\xi) = \max_{k \leq K}  \langle a_k(y), \xi \rangle + b_k$ reduces to a linear program.  One important special case is when the cost function corresponds to the conditional value at risk of a linear function $\langle \bar a(y), \xi \rangle + \bar b$ with confidence level $\alpha \in (0, 1]$ \cite{rockafellar}, where K = 2 written as the piecewise affine function formation $h(y,\xi)$. Then the cost function takes the form
\begin{equation}\label{CVaRDRO}
\begin{split}
J^*  & = \inf_{y \in \mathbb{Y}} \left\{CVaR_\alpha(\langle \bar a(y), \xi \rangle) + \bar b \right\}, \\
 & = \inf_{y \in \mathbb{Y}} \left\{\inf_{\tau \in \mathbb{R}}\mathds{E}^{\mathds{P}}[-\tau + \frac{1}{\alpha} \max \{\langle  \bar a(y),\xi \rangle + \bar b + \tau, 0\} ] \right\},\\
 & = \inf_{y \in \mathbb{Y},\tau \in \mathbb{R}}\mathds{E}^\mathds{P}[\max_{k = 1,2}  \langle a_k(y), \xi \rangle + b_k(\tau)],
\end{split}
\end{equation}
Suppose the support of $\mathds{P}$ is defined by $\Xi :=\{\xi \in \mathbb{R}^m: H\xi \leq d\}$. Then the general form of the distributionally robust optimization problem for CVaR 
\begin{equation}\label{DRO}
\hat{J}_{DRO} = \inf_{y \in \mathbb{Y}, \tau \in \mathbb{R}} \sup_{\mathds{Q} \in \hat{\mathcal{P}}_N}\mathds{E}^\mathds{Q}[\max_{k=1,2} \langle a_k(y),\xi \rangle + b_k],
\end{equation}
as shown in \cite{mohajerin}, can be equivalently reformulated as the linear program
\begin{eqnarray}\label{worstcase}
\inf_{\lambda, s_i, \gamma_{ik}, y, \tau} \quad \lambda\varepsilon + \frac{1}{N_s}\sum_{i=1}^{N_s}s_i,
\end{eqnarray}
\text{subject to}
\begin{eqnarray*}
b_k +  \langle a_k(y), \hat{\xi}_i\rangle + \langle\gamma_{ik}, d-H\hat{\xi}\rangle \le s_i,\ & \forall i\le N_s, \forall k=1,2,\\
 ||H'\gamma_{ik}-a_k(y)||_\infty \le \lambda, & \forall i\le N_s,  \forall k=1,2,\\
 \gamma_{ik} \geq 0,  & \forall i\le N_s, \forall k=1,2,
\end{eqnarray*}
where $\varepsilon \geq 0$. For $\varepsilon \to 0$, the optimization results of \eqref{worstcase} correspond to the expectation of the loss function $h(y, \hat{\xi})$ under the empirical distribution, where $\frac{1}{N_s}\sum_{i=1}^{N_s} s_i$ represents the sample average of the loss function $h(y, \hat{\xi})$.

\subsection{Distributionally Robust Stochastic OPF with CVaR}
We now use the above developments to formulate a data-driven distributionally robust OPF problem to balance an operating efficiency metric with CVaR values of line flow and local device constraint violations. Considering the stochastic OPF shown in \eqref{stochopf}, the $m$-th power line flow constraints \eqref{lineconstraints}  can be written in the form
\begin{equation} \label{linearconstraintgeneral}
g_m(D,e,\xi) = [\Phi(D)]_m\xi + [b(e)]_m,
\end{equation}
where $[.]_m$ denotes the $m$-th elements of a vector or $m$-th row of a matrix. A similar form for the $n$-th local devices constraint can be obtained, which we denote as $h_{jq}(D,e,\xi)$. The decisions variables $D$ and $e$  both appear linearly in $g_m$ and $h_{jq}$. The CVaR of the line flow and local devices constraints are then given by \cite{tyler}
\begin{equation} \label{markovgeneration}
\begin{split}
  & \mathds{E}[g_m(D,e,\xi) + \tau_m ]_+ - \tau_m\alpha \le 0, m = 1,...,2LT\\
  & \mathds{E}[h_{jq}(D,e,\xi)+\tau_{jq}]_+ - \tau_{jq}\alpha \le 0, j = 1,..,N_d, q = 1.,,,. l_j.
\end{split}
\end{equation}
We collect the line flow and local device constraints into a list of $V = 2LT + \sum_{j=1}^{N_d}l_j$ elements. The first $2LT$ elements are line flow constraints, and the remaining $\sum_{j=1}^{N_d}l_j$ elements are local device constraints. The CVaR value of the $v$-th constraint in the above list is denoted by $\Theta_v(D,e,\xi,\tau) = \max_{k=1,2} \langle a_{vk}(y),\xi \rangle + b_{vk}(\tau)$, where decision vector $y$ consists of optimization variables $D$, $e$. 
For $v = 1,...,2LT$, the parameters of the line constraints CVaR $\Theta_v$ are as follows, with $\xi$ redefined here as $[\xi, 1]'$ \\
\begin{equation*}
a_{v1}(y)=\left[
\begin{matrix}
\sum_{j=1}^{N_d}\Gamma_jC_jB_jD_j,  \sum_{j=1}^{N_d}\Gamma_jC_jB_je_j
\end{matrix}\right]
\end{equation*}
\begin{equation*}
b_{v1}(\tau) = -\bar{p} +\sum_{j=1}^{N_d}\Gamma_j(r_j + G_j\xi + C_jA_jx_0^j) + \tau_v - \tau_v\alpha,\\
\end{equation*}
\begin{equation*}
a_{v2}(y)=\left[
\begin{matrix}
0, 0
\end{matrix}\right], 
b_{v2}(\tau) = -\tau_v\alpha
\end{equation*}
For $v = 2LT+1,...,V$, the CVaR of the local device constraints $\Theta_v$ can be defined analogously to obtain similar forms for $a_{v1}(y), b_{v1}(\tau)$ and $a_{v2}(y), b_{v2}(\tau)$.

Based on the above, we now formulate a distributionally robust OPF problem \eqref{DROOPF}, which minimizes a weighted sum of operation cost and the CVaR of the line and local device constraint violations:
\begin{equation} \label{DROOPF}
\begin{split}
& \hat{J}_{OPF} = \hat{J}_{Cost} + \rho\hat{J}_{DRO}\\
& = \hspace{-1mm} \inf_{\begin{subarray}{1}y \in \mathbb{Y} \\ \tau \in \mathbb{R} \end{subarray}} \bigg\{\mathds{E}\sum_{j=1}^{N_d} J_j(\hat{\xi}) +  \rho \hspace{-2 mm}\sup_{\mathds{Q} \in \hat{\mathcal{P}}_N}\sum_{v=1}^V \mathds{E}^\mathds{Q}[\Theta_v(D,e,\xi,\tau)]  \bigg\}, \\
& =\hspace{-1mm} \inf_{\begin{subarray}{1} y \in \mathbb{Y}, \tau \in \mathbb{R}, \\ \lambda_v, s_{iv}, \gamma_{ikv} \end{subarray}} \bigg\{\mathds{E}\sum_{j=1}^{N_d} J_j(\hat{\xi}) + \sum_{v=1}^V \bigg(\lambda_v\varepsilon + \frac{1}{N_s}\sum_{i=1}^{N_s} s_{iv}\bigg) \bigg\},\\
\end{split}
\end{equation}
\text{subject to}:
\begin{eqnarray*}
\sum_{j=1}^{N_d} (r_j + C_j(A_jx_{0}^{j}+B_je_j)) = 0,\\
\sum_{j=1}^{N_d}(G_j + C_jB_jD_j)=0, \\
\rho ( b_{kv}(\tau) + \langle a_{kv}(y), \hat{\xi}_i\rangle + \langle\gamma_{ikv}, d-H\hat{\xi}_i\rangle) \le s_{iv},&\\
 ||H'\gamma_{ikv}-\rho a_{kv}(y)||_\infty \le \lambda_v,&\\
 \gamma_{ikv} \geq 0,&\\
 \forall i\le N_s, \forall v \le V, k=1,2 &\\
\end{eqnarray*}
where $\rho \in \mathbb{R}_+$ quantifies the power system operators' risk aversion. This is a quadratic program that explicitly uses the training dataset $\hat{\Xi}_N := \{\hat{\xi_i}\}_{i\leq N}$. The risk aversion parameter $\rho$ and the Wasserstein radius $\varepsilon$ allow us to explicitly balance tradeoffs between efficiency, risk, and sampling errors inherent in $\hat{\Xi}_N$.

\subsection{Data-Driven Distributionally Robust Model Predictive Control}
The proposed distributionally robust optimization framework can be utilized as a component of a multi-stage closed-loop model predictive control (MPC) strategy. This strategy allows explicit incorporation of historical datasets and forecasting techniques that can update predictions during operation. MPC is an feedback control technique that implements optimal open-loop control decisions at the current time step $t$-th while considering future system states over a fixed time horizon $\mathcal{H}_t$ \cite{camacho2013model,fortenbacher2014grid,adamek2014decisive,parisio2016stochastic}. In the above present problem, the distributionally robust stochastic OPF \eqref{DROOPF} could be regarded as a building block for a MPC-based closed-loop control strategy, which includes the current and future system states based on the forecasts of the renewable energy availability and the local devices dynamics. The MPC-based control strategy involves the following steps:
\begin{itemize}
\item At time step $t$, forecast the uncertainties across the power network over horizon $\mathcal{H}_t$.
\item Solve \eqref{DROOPF} over horizon $\mathcal{H}_t$.
\item Implement the reserve control policies for each device at current time step $t$.
\item Move to time step $t+1$, and return to the first step.
\end{itemize}
The optimal reserve control policies are computed for the entire time horizon $\mathcal{H}_t$, but only the decisions for the current time step are implemented. Then we shift to next time step with the forecasts updated. To simplify the exposition, our numerical results will focus on solving \eqref{DROOPF} in one particular time step. We are pursuing and elaborating on a multi-stage data-driven distributionally robust MPC algorithm and properties in ongoing work.

\section{Numerical Example}
To illustrate the proposed framework, we use a modified IEEE 118 bus test network as defined in \cite{118testbed}  shown in Fig.1. A wind farm is connected to bus 9 with nominal infeed of 1000 MW. To understand the trade offs of the proposed distributionally robust optimization methodology, we consider only the line flow constraint between bus 8 and bus 9, marked by a cross in Fig.1. The maximum rating on the power line between bus 8 and 9 is 950 MW. Other line flow constraints, security constraints, and voltage constraints can also be easily incorporated, and we are pursuing this in ongoing work. The wind power forecast errors are derived from real wind data from the hourly wind power measurements provided in the 2012 Global Energy Forecasting competition (GEFCom2012) \cite{wind}. The wind power forecast errors are based on the so-called persistence forecast, which predicts the wind power output at the next time step to be equal to that at the previous time step. It can be seen that the forecast errors are highly leptokurtic, i.e., that the errors are have outliers that make the distribution tails much heavier than Gaussian tails. Since the wind power data from GEFCom2012 are normalized in [0,1], we scale the forecast errors to have zero mean and the standard deviation $\sigma = 300$ MW.
In this case study, we only consider a single line constraint, hence no other local device constraint is included. Additionally, we do not assume a bound on the wind power forecast errors $\xi$ , therefore, the $H$ and $d$ in \eqref{DROOPF} are equal to zero.
\begin{figure}[h]
\centering
\includegraphics[width=3.5in]{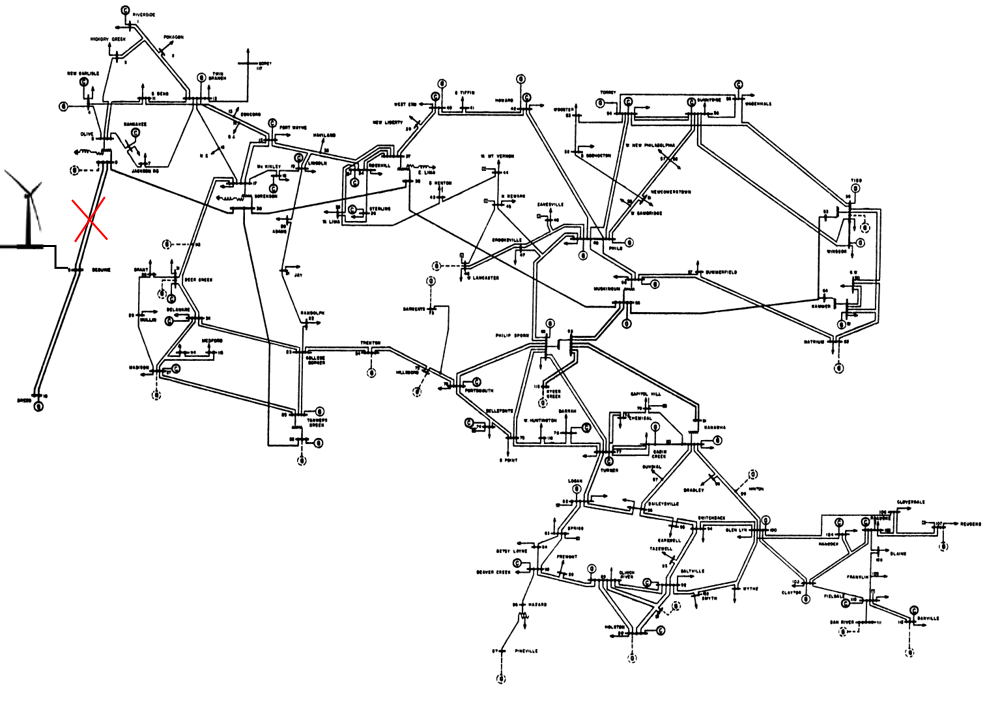}
\caption{IEEE 118 bus test network.}
\label{fidg:side:b}
\end{figure}

For comparison purposes, we also introduce two other methods which we call CVaR OPF and Gaussian OPF to evaluate the constraint violations. In CVaR OPF, we use a basic sample average to approximate CVaR of the line constraint violations with training dataset $\hat{\Xi}_{N_s}$, which corresponds to setting $\varepsilon=0$ in \eqref{DROOPF}. In Gaussian OPF, we assume the data are generated by a Gaussian distribution with the sample mean and variance from dataset $\hat{\Xi}_{N_s}$. We subsample $N_s=100$ elements from the training dataset $\hat{\Xi}_{N_s}$.
\begin{figure}[h]
\centering
\includegraphics[trim=00mm 0mm 00mm 0mm, clip, width=0.48\textwidth]{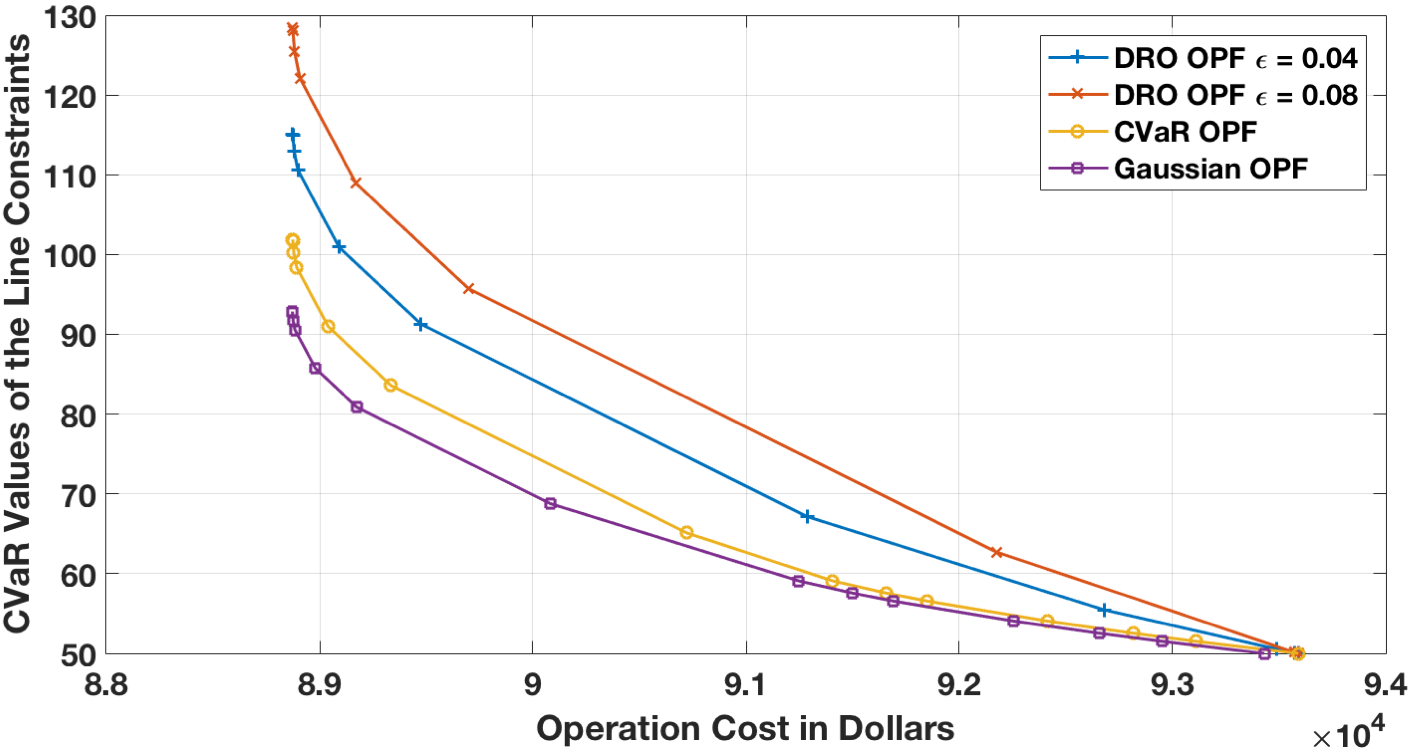}
\caption{Predicted tradeoff between operating cost and CVaR of line constrain violations, based on different models of data uncertainty. The Gaussian model underestimates the risk, as demonstrated by the out-of-sample evaluations in Fig.3.}
\label{fidg:side:b}
\end{figure}

\begin{figure}[h]
\centering
\includegraphics[trim=10mm 0mm 16mm 0mm, clip, width=0.48\textwidth]{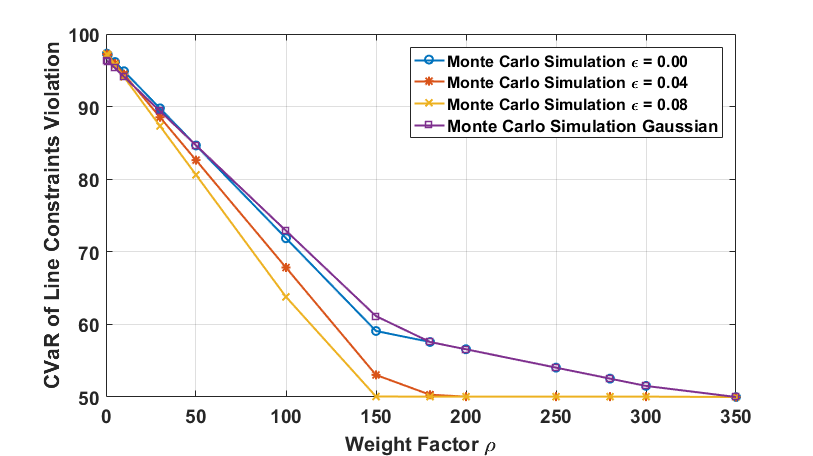}
\caption{Out-of-sample performance is superior to other approaches and conservativeness is explicitly controllable, as demonstrated by Monte Carlo simulations. For each value of $\rho$ we subsampled new values from the training dataset and computed the empirical CVaR of the line constraint.}
\label{fidg:side:b}
\end{figure}

\begin{figure}[h]
\centering
\includegraphics[trim=10mm 0mm 15mm 0mm, clip, width=0.48\textwidth]{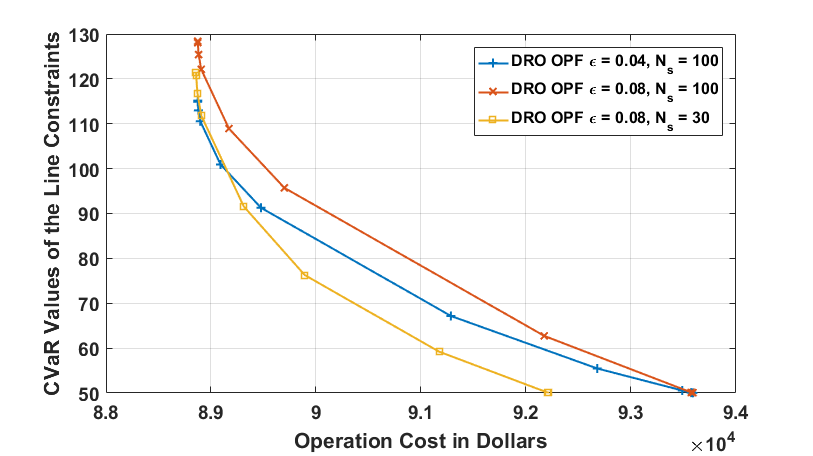}
\caption{Sample errors can lead to the underestimation of the risk of constraint violations when the data contain less variance than the real distribution, but their effect is diminished by the data-driven distributionally robust approach.}
\label{fidg:side:b}
\end{figure}

Fig.2 visualizes the solutions of the proposed stochastic OPF by three methods: CVaR OPF, Gaussian OPF and DRO OPF ($\varepsilon = 0.04, 0.08$). The results illustrate the tradeoffs between optimal operation cost and CVaR of line constraint violations by varying the risk aversion parameter $\rho$. The tradeoffs could be explicitly adjusted by $\rho$ according to the system operators' efficiency and risk tradeoff preferences.

Fig.2 also demonstrates the underestimation of the risk of constraint violations by assuming a Gaussian distribution. Distributionally robust optimization considers the worst-case scenario of the uncertainties to produce more conservative decisions on generator power adjustments and reserve policies. The conservativeness of the results are controllable by adjusting the Wasserstein distance $\varepsilon$. In this paper, we presents two cases $\varepsilon = 0.04, \varepsilon = 0.08$. To emphasize, the results of DRO OPF ($\varepsilon = 0.00$) overlap with the curve given by CVaR OPF in Fig.2. Increasing $\varepsilon$ provides better robustness to sampling errors and yields superior out-of-sample performance guarantees.

Fig.3 visualizes the out-of-sample performance guarantee of stochastic OPF solution offered by distributionally robust model based on Monte Carlo simulation. For each value of $\rho$ we subsampled new values from the training dataset and computed the empirical CVaR of the line constraint. The decisions from DRO OPF method ensure smaller line constraint violation for all values of the risk aversion $\rho$. Again, we also see that increasing the Wasserstein radius $\varepsilon$ provides lower risk. The benefits saturate for small $\rho$.

In the DRO OPF case, we also re-solved stochastic OPF problem with reduced size of training dataset to $N_s = 30$. Fig.4 visualizes the effects of sampling errors. We see in this case that the smaller dataset has caused an underestimation of risk, since the reduced data contains less variance than the full data. Sampling errors can cause either underestimation or overestimation of actual risks depending on whether the dataset features more or less variation than the real data generating distribution, but their effect is diminished with the proposed distributionally robust approach.

\section{Conclusions}
In this paper, we have proposed a distributionally robust optimization method to solve a stochastic OPF problem with the limited knowledge of the uncertain forecast error across the power system network. By explicitly quantifying the differences between the empirical data-based distributions and real data-generating distributions using the Wasserstein metric, we computed power schedules and reserve policies which were robust to sampling errors. This approach has superior out-of-sample performance guarantees in comparison with other approaches and allows a more systematic trade off of efficiency and risk. A simple numerical example illustrated these basic tradeoffs. The results also demonstrated that making inaccurate assumptions about probability distributions can underestimate the risk of constraint violations compared to the distributionally robust optimization. Ongoing and future work includes case studies in larger networks and adapting the approach to handle voltage constraints in distribution networks using alternative linearizations of the power flow equations, as in \cite{DallAnese2}.
\bibliography{GuoBakerDallaneseHuSummers}
\bibliographystyle{ieeetr}

\end{document}